\documentclass[12pt]{article}
\usepackage{amsmath,amssymb,amsbsy,amsfonts,amsthm,latexsym,
            amsopn,amstext,amsxtra,euscript,amscd}

\newtheorem{lem}{Lemma}
\newtheorem{lemma}[lem]{Lemma}

\newtheorem{thm}{Theorem}
\newtheorem{theorem}[thm]{Theorem}

\def\\{\cr}
\def\({\left(}
\def\){\right)}
\def\[{\left[}
\def\]{\right]}
\def\<{\langle}
\def\>{\rangle}

\begin{document}

\title{Counting Keith numbers}

\author{
{\sc Martin~Klazar} \\
{Department of Applied Mathematics and}\\
{Institute for Theoretical Computer Science (ITI)}\\
{Faculty of Mathematics and Physics, Charles University} \\
{Malostransk\'e n\'am. 25, 11800 Praha, Czech Republic} \\
{\tt klazar@kam.mff.cuni.cz}\\
and\\
{\sc Florian~Luca} \\
{Instituto de Matem{\'a}ticas}\\
{Universidad Nacional Autonoma de M{\'e}xico} \\
{C.P. 58089, Morelia, Michoac{\'a}n, M{\'e}xico} \\
{\tt fluca@matmor.unam.mx}
}

\date{\today}
\maketitle

\begin{abstract}
A Keith number is a positive integer $N$ with the decimal 
representation $a_1a_2\dots a_n$ 
such that $n\ge 2$ and $N$ appears in the sequence $(K_m)_{m\ge 1}$ 
given by the 
recurrence 
$K_1=a_1,\dots,K_n=a_n$ and 
$K_m=K_{m-1}+K_{m-2}+\cdots+K_{m-n}$ for $m>n$.
We prove that there are only finitely many Keith numbers using only one 
decimal digit (i.e., $a_1=a_2=\cdots=a_n$),
 and that the set of Keith numbers is of asymptotic density zero.
\end{abstract}

\section{Introduction}

With the number $197$, let $(K_m)_{m\ge 1}$ be the sequence 
whose first three terms $K_1=1,~K_2=9$ and $K_3=7$ are the digits of $197$ 
and which satisfies the recurrence $K_m=K_{m-1}+K_{m-2}+K_{m-3}$ for all 
$m>3$. Its initial terms are 
$$
1,~9,~7,~17,~33,~57,~107,~197,~361,~665,\ldots
$$
Note that $197$ itself is a member of this sequence. This phenomenon 
was first noticed 
by Mike Keith and such numbers are now called {\it Keith numbers}.
More precisely, a number $N$ with decimal representation
$a_1a_2\ldots a_n$ is a Keith number if $n\ge 2$ and $N$  
appears in  the sequence $K^N=(K^N_m)_{m\ge 1}$ whose $n$ initial terms 
are the digits of $N$ read from left to right 
and satisfying $K^N_m=K^N_{m-1}+K^N_{m-2}+\cdots+K^N_{m-n}$ 
for all $m>n$.  These numbers appear in Keith's papers \cite{Ke1} and 
\cite{Ke2} and they are the subject of entry $A007629$ in 
Neil Sloane's Encyclopedia of Integer Sequences \cite{Slo} (see also 
\cite{Pick1}, \cite{Pick2} and \cite{Pick3}). 

\medskip

Let ${\cal K}$ be the set of all Keith numbers. It is not known  if 
${\cal K}$ is infinite or not. The 
sequence ${\cal K}$ begins
$$
14,~19,~28,~47,~61,~75,~197,~742,~1104,~1537,~2208,
2580,~3684,~4788,\dots
$$
In total there are $94$ Keith numbers smaller than $10^{29}$ (\cite{Ke2}).  
Recall that a rep-digit is a positive 
integer $N$ of the form $a(10^n-1)/9$ for some $a\in \{1,\ldots,9\}$ 
and $n\ge 1$; i.e., 
a number which is a string of the same digit $a$ when written in base 
$10$. Our first result shows that there are only finitely many Keith
numbers which are rep-digits.

\begin{theorem}
\label{thm:krepu}
There are only finitely many Keith
numbers which are rep-digits and their set can be effectively determined.
\end{theorem}

We point out that some authors refer to the Keith numbers as 
{\it replicating Fibonacci digits} in analogy with the Fibonacci sequence 
$(F_n)_{n\ge 1}$ given by $F_1=1,~F_2=1$ and $F_{n+2}=F_{n+1}+F_n$ 
for all $n\ge 1$. In \cite{Lu2} it is shown that the largest
rep-digit Fibonacci number is $55$.  

\medskip

The proof of Theorem \ref{thm:krepu} uses Baker's type estimates
for linear forms in logarithms. It will be clear 
from the proof that it applies to all 
{\it base $b$ Keith numbers} for any fixed integer $b\ge 3$, where these 
numbers are defined analogously starting with their base $b$ expansion 
(see the remark after the proof of Theorem \ref{thm:krepu}). 

\medskip

For a positive integer $x$ we write ${\cal K}(x)={\cal K}\cap [1,x]$. 
As we mentioned before, ${\cal K}(10^{29})=94$.
A heuristic argument in \cite{Ke2} suggests that 
$\#{\cal K}(x)\gg \log x$, and, in particular, that ${\cal K}$ should be 
infinite. Going in the opposite way, we show that ${\cal K}$ is of asymptotic 
density zero.

\begin{theorem}
\label{thm:countKeith}
The estimate 
$$
\#{\cal K}(x)\ll \frac{x}{\sqrt{\log x}}
$$
holds for all positive integers $x\ge 2$.
\end{theorem}

The above estimate is very weak. It does not even imply that 
that sum of the reciprocals of the members of ${\cal K}$ is convergent.
We leave to the reader the task of finding a better upper bound on 
$\#{\cal K}(x)$. Typographical changes (see the remark after 
the proof of Theorem \ref{thm:countKeith}) show that Theorem 
\ref{thm:countKeith} also is valid for the set 
of base $b$ Keith numbers if $b\ge 4$. Perhaps 
it can be extended also to the case
$b=3$. For $b=2$, Kenneth Fan has an unpublished 
manuscript showing how to construct all Keith numbers 
(see \cite{Ke2}) and that,
in particular, there are infinitely many of them. For example, any power of 
$2$ is a binary Keith number.

\medskip

Throughout this paper, we use the Vinogradov symbols $\gg$ and $\ll$ as well
as the Landau symbols $O$ and $o$ with their usual meaning. 
Recall that for functions $A$ and $B$ the inequalities 
$A\ll B$, $B\gg A$ and $A=O(B)$ are all equivalent to the 
fact that there exists a positive constant $c$ such that the inequality 
$|A|\le cB$ holds. The constants in the inequalities implied by these 
symbols 
may occasionally depend on other parameters. For a real number $x$ we 
use $\log x$ for the natural logarithm of $x$. For a set ${\cal A}$, we use 
$\#{\cal A}$ and $|{\cal A}|$ to denote its cardinality.

\medskip

{\bf Acknowledgements.} The first author acknowledges gratefully the 
institutional 
support to ITI by the grant 1M0021620808  of
the Czech Ministry of Education. The second author was working on this paper
during a visit to CRM in Montreal during Spring 2006. 
The hospitality and support of
this institution is gratefully acknowledged. During the
preparation of this paper, he was  also supported in
part by Grants SEP-CONACyT 46755, PAPIIT IN104005 and a Guggenheim
Fellowship.

\section{Preliminary Results}

For an integer $N>0$, recall the definition of the sequence 
$K^N=(K^N_m)_{m\ge 1}$ given in the Introduction. In $K^N$ 
we allow $N$ to be 
any string of the digits $0,1,\dots,9$, so $N$ may have initial zeros. 
So, 
for example, $K^{020}=(0,2,0,2,4,6,12,22,\dots)$. For $n\ge 1$ we define the 
sequence $L^n$ 
as $L^n=K^M$ where $M=11\dots 1$ with $n$ digits $1$. In particular, 
$L^1=(1,1,1,\dots)$ and $L^2=(1,1,2,3,5,8,\dots)$, the Fibonacci numbers. 
In the following lemma, which will be used in the proofs of both Theorems 
1 and 2, 
we establish some properties of the sequences $K^N$ and $L^n$.

\begin{lemma}\label{simplelemma}
Let $N$ be a string of the digits  $0,1,\dots,9$ with length $n\ge 1$. If 
$N$ 
does not start with $0$, we understand it also as the decimal 
representation of a
positive integer. 

\begin{enumerate} 
\item If $N$ has at least $k\ge 1$ nonzero entries, then 
$K^N_m\ge L^k_{k+m-n}$ holds for every $m\ge n+1$.
\item If $N$ has at least one nonzero entry, then $K^N_m\ge L^n_{m-n}$ holds 
for every $m\ge n+1$. We have $K^N_m\le 9L^n_m$ for every $m\ge 1$.
\item If $n\ge 3$ and $N=K^N_m$ for some $m\ge 1$ (so $N$ is a Keith number), 
then 
$2n<m<7n$.
\item For fixed $n\ge 2$ and growing $m\ge n+1$,
$$
L^n_m=2^{m-n-1}(n-1)(1+O(m/2^n))+1
$$
where the constant in $O$ is absolute.
\end{enumerate}
\end{lemma}
\begin{proof}
1. By the recurrences defining $K^N$ and $L^k$, 
the inequality clearly holds for 
the first $k$ indices $m=n+1,n+2,\dots,n+k$. For $m>n+k$ it holds by 
induction.

\medskip

2. We have $K^N_m\ge 1=L^n_{m-n}$ for $m=n+1,n+2,\dots,2n$ and the inequality 
holds. For $m>2n$ it holds by induction. The second inequality follows 
easily by induction.

\medskip

3. The lower bound $m>2n$ follows from the fact that $K^N$ is nondecreasing 
and that 
$$
K^N_{2n}\le 9L^n_{2n}=9\cdot2^{n-1}(n-1)+9<10^{n-1}\le N
$$ 
for $n\ge 3$. To obtain the upper bound, 
note that for $m\ge n$ we have by induction that 
$L^n_m\ge L^2_{m-n+2}\ge\phi^{m-n}$ where $\phi=1.61803\dots$ 
is the golden ratio. Thus, by part 2, 
$$
10^n>N=K^N_m\ge L^n_{m-n}\ge\phi^{m-2n}
$$
and $m<(2+\log 10/\log\phi)n<7n$.

\medskip

4. We write $L^n_m$ in the form $L^n_m=(2^{m-n-1}-d(m))(n-1)+1$ and prove 
by induction on $m$ that for $m\ge n+1$, 
$$
0\le d(m)<m2^{m-2n}.
$$
This will prove the claim.

\medskip

It is easy to see by the recurrence that 
$L^n_{n+1},L^n_{n+2},\dots, L^n_{2n+1}$ are equal, respectively, 
to $2^0(n-1)+1,2^1(n-1)+1,\dots,2^n(n-1)+1$. So $d(m)=0$ for 
$n+1\le m\le 2n+1$ and the claim holds. 
For $m\ge 2n+1$,
\begin{eqnarray*}
L^n_m&=&L^n_{m-1}+L^n_{m-2}+\cdots+L^n_{m-n}\\
&=&\sum_{k=1}^n\Big((2^{m-n-1-k}-d(m-k))(n-1)+1\Big)\\
&=&\Big(2^{m-n-1}-2^{m-2n-1}+1-\sum_{k=1}^n d(m-k)\Big)(n-1)+1
\end{eqnarray*}
and the induction hypothesis give
\begin{eqnarray*}
0\le d(m)&=&2^{m-2n-1}-1+\sum_{k=1}^n d(m-k)\\
&<&2^{m-2n-1}+(m-1)\sum_{k=1}^n 2^{m-2n-k}\\
&<&m2^{m-2n}.
\end{eqnarray*}
\end{proof}

\noindent

In part 4, if $m$ is roughly of size $2^n$ and 
larger then the error term swallows 
the main term and the asymptotics is useless. Indeed, the correct 
asymptotics of $L^n_m$ 
when $m\to\infty$ is $c\alpha^m$ where $c>0$ is a constant and $\alpha<2$ 
is the only positive root of the polynomial $x^n-x^{n-1}-\cdots-x-1$. 
But for $m$ small relative to $2^n$, say 
$m=O(n)$ (ensured for Keith numbers by part 3), this ``incorrect'' 
asymptotics of $L^n_m$ is very precise and useful, as we shall demonstrate in 
the proofs of Theorems \ref{thm:krepu} and \ref{thm:countKeith}.

\medskip

In the proof of Theorem \ref{thm:krepu} we will apply also a 
lower bound for a linear form in logarithms.
The following result can be deduced from Corollary 2.3 of \cite{Mat}.

\begin{lemma}
\label{lem:linearform}
Let $A_1,\dots,A_k$, $A_i>1$, and $n_1,\dots,n_k$ 
be integers, and let $N=\max\{|n_1|,\ldots,|n_k|,2\}$. 
There exist positive absolute constants $c_1$ and $c_2$ (which are 
effective), such that if
$$
\Lambda=n_1\log A_1+n_2\log A_2+\cdots+n_k\log A_k\ne 0,
$$
then
$$
\log |\Lambda|>-c_1 c_2^k (\log A_1)\ldots (\log A_k)\log N.
$$
\end{lemma}

For the proof of Theorem 2 we will need an upper bound on sizes of antichains
(sets of mutually incomparable elements) in the poset (partially ordered set)
$$
P(k,n)=(\{1,2,\dots,k\}^n,\le_p)
$$ 
where $\le_p$ is the product ordering
$$
a=(a_1,a_2,\dots, a_n)\le_p b=(b_1,b_2,\dots, b_n)\iff a_i\le 
b_i\mbox{ for } i=1,2,\dots,n.
$$
We have $|P(k,n)|=k^n$ and for $k=2$ the poset $P(2,n)$ is the Boolean poset 
of subsets of an $n$-element set ordered by inclusion. 
The classical theorem of 
Sperner (see \cite{aign_zieg} or \cite{sperner_theory}) 
asserts that the maximum size of 
an antichain in $P(2,n)$ equals to 
the middle binomial coefficient 
$\binom{n}{\lfloor n/2\rfloor}$.
In the next lemma we obtain an upper bound for any $k\ge 2$.

\begin{lemma}\label{antichains}
If $k\ge 2, n\ge 1$ and $X\subset P(k,n)$ is an antichain to $\le_p$, then 
$$
|X|<\frac{(k/2)\cdot k^n}{n^{1/2}}.
$$
\end{lemma} 
\begin{proof}
We proceed by induction on $k$. 
For $k=2$ this bound holds by Sperner's theorem 
because
$$
\binom{n}{\lfloor n/2\rfloor}<\frac{2^n}{n^{1/2}} 
$$
for every $n\ge 1$. Let $k\ge 3$ and $X\subset P(k,n)$ be an antichain. 
For $A$ running through the subsets of $[n]=\{1,2,\dots,n\}$, we partition 
$X$ in the sets $X_A$ where $X_A$ consists of the $u\in X$ satisfying 
$u_i=k\iff i\in A$.
If we delete from all $u\in X_A$ all appearances of $k$, we obtain 
(after appropriate relabelling of coordinates) a set of $|X_A|$ distinct 
$(n-|A|)$-tuples from $P(k-1,n-|A|)$ that must be an antichain to $\le_p$. 
Thus, by induction, for $|A|<n$ we have
$$
|X_A|<\frac{((k-1)/2)\cdot(k-1)^{n-|A|}}{(n-|A|)^{1/2}}
$$
and $|X_{[n]}|\le 1$.
Summing over all $A$s and using the inequality $\sqrt{n/m}\le(n+1)/(m+1)$ 
(which holds for $1\le m\le n$) and standard properties of binomial 
coefficients, we get
\begin{eqnarray*}
|X|&=&\sum_{A\subset[n]}|X_A|\\
&<&1+\sum_{i=0}^{n-1}\binom{n}{i}\frac{((k-1)/2)
\cdot(k-1)^{n-i}}{(n-i)^{1/2}}\\
&=&\frac{1}{\sqrt{n}}\left(\sqrt{n}+\frac{1}{2}
\sum_{i=0}^{n-1}\binom{n}{i}\sqrt{n/(n-i)}\cdot(k-1)^{n-i+1}\right)\\
&\le&\frac{1}{\sqrt{n}}\left(\sqrt{n}+\frac{1}{2}
\sum_{i=0}^{n-1}\binom{n+1}{n-i+1}(k-1)^{n-i+1}\right)\\
&<&\frac{k^{n+1}}{2\sqrt{n}}.
\end{eqnarray*}
\end{proof}

\noindent
We conclude this section with three remarks as to the last lemma. 

\medskip

1. Various generalizations and strengthenings of Sperner's theorem were 
intensively studied, see, e.g., the book of Engel and Gronau 
\cite{sperner_theory}. Therefore, we do not expect much 
originality in our bound. 

\medskip

2. It is clear that for $k=2$ the exponent $1/2$ of $n$ in the bound of 
Lemma~\ref{antichains} cannot be  increased. The same is true for any 
$k\ge 3$. We briefly sketch a construction of a large antichain when 
$k=3$; for $k>3$ similar constructions can be given. For 
$k=3$ and $n=3m\ge 3$ consider the set $X\subset P(3,n)$ consisting of all 
$u$ which have $i$ $1$s, $n-2i$ $2$s and $i$ $3$s, where $i=1,2,\dots,m=n/3$. 
It follows that $X$ is an antichain and that
$$
|X|=\sum_{i=1}^m\binom{n}{i,i,n-2i}=\sum_{i=1}^m\frac{n!}{(i!)^2(n-2i)!}.
$$
By the usual estimates of factorials, if $m-\sqrt{n}<i\le m$ then
$$
\binom{n}{i,i,n-2i}\gg \binom{n}{m,m,m}\gg\frac{3^n}{n}.
$$
Hence $X$ is an antichain in $P(3,n)$ with size
$$
|X|\gg\sqrt{n}\cdot\frac{3^n}{n}=\frac{3^n}{\sqrt{n}}.
$$

\medskip

3. For composite $k$ we can decrease the factor $k/2$ in the bound of 
Lemma~\ref{antichains}. Suppose that $k=lm$ where $l\ge m\ge 2$ are integers 
and let 
$X\subset P(k,n)$ be an antichain. We associate with every $u\in X$ the 
pair of $n$-tuples $(v^u,w^u)\in P(m,n)\times P(l,n)$ defined by 
$v^u_i=u_i-m\lceil u_i/m\rceil+m$ and $w^u_i=\lceil u_i/m\rceil$, 
$1\le i\le n$. Note that the pair $(v^u,w^u)$ uniquely determines $u$ and 
that if $w^u=w^{u'}$ then $v^u$ and $v^{u'}$ are incomparable by $\le_p$. 
Thus, by Lemma~\ref{antichains}, for fixed $w\in P(l,n)$ there are less than 
$(m/2)m^n/\sqrt{n}$ elements $u\in X$ with $w^u=w$. The number of $w$s is at 
most $|P(l,n)|=l^n$. Hence
$$
|X|<\frac{(m/2)\cdot m^n}{n^{1/2}}\cdot l^n=\frac{(m/2)\cdot k^n}{n^{1/2}}.
$$      
In particular, if $k$ is a power of $2$ then $|X|<k^n/\sqrt{n}$ for every 
antichain $X\subset P(k,n)$. 

\section{The proof of Theorem \ref{thm:krepu}}
\label{sec:proofofkrepu}

\noindent
Let $N=a(10^n-1)/9=aa\dots a$, $1\le a\le 9$, be a rep-digit. Since $K^N=aL^n$,
$N$ is a Keith number if and only if the repunit $M=(10^n-1)/9=11\dots 1$ 
is a Keith number. Suppose that $M$ is a Keith number: for some $m$ we have 
$$
M=\frac{10^n-1}{9}=L^n_m=2^{m-n-1}(n-1)\left(1+O\left(\frac{m}{2^n}\right)
\right),
$$
where the asymptotics was proved in Lemma~\ref{simplelemma}.4.
We rewrite this relation as 
$$
\frac{2^{2n+1-m}5^{n}}{9(n-1)}-1=\frac{1}{9(n-1)2^{m-n-1}}+
O\left(\frac{m}{2^n}\right).
$$
Since $2n<m<7n$ by Lemma~\ref{simplelemma}.3, we get
$$
\frac{2^{2n+1-m}5^n}{9(n-1)}-1=O\left(\frac{n}{2^n}\right).
$$
Because $5^n>9(n-1)$ for every $n\ge 1$, the left side is always non-zero. 
Writing it in the form 
$e^{\Lambda}-1$ and using that $e^{\Lambda}-1=O(\Lambda)$ (as $\Lambda\to 0$), 
we get
$$
0\ne\Lambda=(2n+1-m)\log 2+n\log 5-\log(9(n-1))\ll\frac{n}{2^n}.
$$
Taking logarithms and applying Lemma \ref{lem:linearform}, we finally obtain 
$$
-d(\log n)^2<\log|\Lambda|<c(\log n-n\log 2)
$$
where $c,d>0$ are effectively computable constants. 
This implies that $n$ is effectively bounded and completes the proof of 
Theorem~\ref{thm:krepu}.
\qed

\medskip

\noindent {\bf Remark.} The same 
argument shows that for every integer $b\ge 3$ there are only effectively 
finitely many base $b$ rep-digits, i.e., 
positive integers of the form $a(b^n-1)/(b-1)$ with $a\in \{1,\ldots,b-1\}$,
which are base $b$ Keith numbers. Indeed, we argue as for $b=10$ and derive 
the equation
$$
\frac{b^n}{(b-1)(n-1)2^{m-n-1}}-1=O(n/2^n).
$$
In order to apply Lemma \ref{lem:linearform}, we need to justify that the 
left side is not zero. If $b$ is not a power of $2$, it has an odd prime 
divisor $p$, and $p^n$ cannot be cancelled, for big enough $n$, by 
$(b-1)(n-1)$. If $b\ge 3$ is a power of $2$, then $b-1$ is odd and has an 
odd prime divisor, which cannot be cancelled by the rest of the expression.

\section{The proof of Theorem \ref{thm:countKeith}}

\noindent
For an integer $N>0$, we denote by $n$ the number of its digits: 
$10^{n-1}\le N<10^{n}$. We shall prove that there are $\ll 10^n/\sqrt{n}$ 
Keith numbers with $n$ digits; it is easy to see that this implies Theorem 2.
There are only few numbers with $n$ digits and $\ge n/2$ zero digits: their
number is bounded by
$$
\sum_{i\ge n/2}\binom{n}{i}9^{n-i}\le n2^n9^{n/2}=n6^n\ll (10^n)^{0.8}.
$$
Hence it suffices to count only the Keith numbers with $n$ digits, 
of which at least half are nonzero.

\medskip

Let $N$ be a Keith number with $n\ge 3$ digits, at least half of them nonzero. 
So, $N=K^N_m$ for some index $m\ge 1$. By Lemma~\ref{simplelemma}.3, 
$2n<m<7n$ and we may use the asymptotics in Lemma~\ref{simplelemma}.4.
Setting $k=\lfloor n/2\rfloor$ and using the inequality in 
Lemma~\ref{simplelemma}.1, we get 
$$
10^n>N=K^N_m\ge L^k_{k+m-n}.
$$
Lemma~\ref{simplelemma}.4 gives that for big $n$,
$$
L^k_{k+m-n}>\frac{2^{m-n-1}(k-1)}{2}>\frac{2^{m-n}n}{12}.
$$ 
On the other hand, the second inequality in Lemma~\ref{simplelemma}.2 
and Lemma~\ref{simplelemma}.4 give, for big $n$, 
$$
10^{n-1}\le N=K^N_m\le 9L^n_m<9\cdot 2^{m-n}n.
$$
Combining the previous inequalities, we get
$$
\frac{10^n}{90}<2^{m-n}n<12\cdot 10^n.
$$
This implies that, for $n>n_0$, the index $m$ attains at most $12$ 
distinct values and
$$
m=(1+\log 10/\log 2+o(1))n=(\kappa+o(1))n.
$$   

\medskip

Now we partition the set $S$ of considered Keith numbers 
(with $n$ digits, at least half of them nonzero) in blocks of numbers $N$ 
having the same value of the index $m$ and the same string of the first 
(most significant) $k=\lfloor n/2\rfloor$ digits.  So, 
we have at most $12\cdot 10 ^k$ blocks. We show in a moment that the numbers 
in one block $B$, when regarded as $(n-k)$-tuples from $P(10,n-k)$, 
form an antichain to $\le_p$. Assuming this, 
Lemma~\ref{antichains} implies that $|B|<10^{n-k+1}/2\sqrt{n-k}$.
Summing over all blocks, we get
$$
|S|<12\cdot 10 ^k\cdot \frac{10^{n-k+1}}{2\sqrt{n-k}}\ll
\frac{10^{n}}{\sqrt{n}},
$$
which proves Theorem 2.

\medskip

To show that $B$ is an antichain, we suppose for the contradiction that $N_1$ 
and $N_2$ are two Keith numbers from $B$ with $N_1<_p N_2$. Let $M=N_2-N_1$
and $M^*=00\dots 0M\in P(10,n)$ (we complete $M$ to a string of length $n$ by 
adding initial zeros). It follows that $M$ has at most $n-k$ digits and 
$M<10^{n-k}$. On the other hand, by the linearity of recurrence and by 
$N_1<_p N_2$, we have 
$$
M=N_2-N_1=K^{N_2}_m-K^{N_1}_m=K^{M^*}_m.
$$
Since $M^*$ has some nonzero entry, the first inequality in 
Lemma~\ref{simplelemma}.2 and Lemma~\ref{simplelemma}.4 give,
for big $n$,
$$
K^{M^*}_m\ge L^n_{m-n}>2^{m-2n-2}n.
$$
Thus 
$$
10^{n-k}=10^{n-\lfloor n/2\rfloor}>M>2^{m-2n-2}n.
$$
Using the above asymptotics of $m$ in terms of $n$, we arrive at the 
inequality
\begin{eqnarray*}
\exp(({\textstyle\frac{1}{2}}\log 10+o(1))n)&>&
\exp((\kappa\log 2-2\log2+o(1))n)\\
&=&\exp((\log 5+o(1))n)
\end{eqnarray*}
that is contradictory for big $n$ because $10^{1/2}<5=10/2$. 
This finishes the proof of Theorem 2.
\qed 

\medskip

\noindent {\bf Remark.} 
The above proof generalizes, with small
modifications, to all bases $b\ge 4$.
We replace base $10$ by $b$, modify the proof accordingly, and have to 
satisfy two conditions. First,  in the beginning of the proof we delete 
from the numbers with $n$ base $b$ digits those with 
$>\alpha n$ zero digits, for some constant $0<\alpha<1$. In order that we 
delete negligibly many, compared to $b^n$, numbers, we must have 
$2\cdot (b-1)^{1-\alpha}<b$. Second, for the final contradiction 
we need that $b^{\alpha}<b/2$. For $b\ge 5$, both conditions are satisfied 
with $\alpha=1/2$, as in case $b=10$. For $b=4$ they are satisfied with 
$\alpha=0.49$, say. However, for $b=3$ they cannot be satisfied by any 
$\alpha$. Thus, the case $b=3$ seems 
to require more substantial changes.

\end{document}